\documentclass[11pt,a4paper]{article} 

\usepackage{anysize}
\usepackage[format = hang]{caption}
\usepackage{amsmath,amsthm,verbatim,amssymb,amsfonts,amscd, graphicx}
\usepackage{graphics,natbib}
\usepackage{titlesec,footmisc}
\usepackage{hyperref} 

\usepackage{listings}
\usepackage{booktabs} 

\theoremstyle{plain}

\theoremstyle{definition}

\newcommand{\cov}{\text{cov}}

\def\Var{{\rm Var}}

\def\pr{{\rm pr}}

\def\hatt{\widehat}
\def\arr{\rightarrow}
\def\N{{\rm N}}
\def\half{\hbox{$1\over 2$}}
\def\eps{\varepsilon}

\def\beq{\begin{eqnarray}}
\def\eeq{\end{eqnarray}}

\def\quantile{{\rm quant}}

\def\beqn{\begin{eqnarray*}}  
\def\eeqn{\end{eqnarray*}}

\def\Var{{\rm Var}}

\def\N{{\rm N}}

\def\Pr{P}
\def\pr{{\rm pr}}

\def\arr{\rightarrow}
\def\hatt{\widehat}

\def\eps{\varepsilon}
\def\half{\hbox{$1\over2$}}

\def\quart{\hbox{$1\over4$}}

\def\rootn{\sqrt{n}}

\def\med{{\rm med}}

\titleformat{\section}{\normalfont\large\sc\centering}{\thesection}{1em}{}
\titleformat{\subsection}[runin]{\normalfont\large\bfseries}{\thesubsection}{1em}{}
\numberwithin{equation}{section} 
\renewenvironment{abstract}
               {\list{}{\rightmargin\leftmargin}%
                \item[\text{\hspace{10mm}\sc Abstract.}]\relax}
               {\endlist}

\begin{document}

\def\heute{November 2026}

\begingroup
\begin{centering} 
\Large{\bf Variance estimators
    robust against jumps and fluctuations}\\[0.8em]
\large{\bf C\'eline Cunen$^1$ and Nils Lid Hjort$^2$} \\[0.3em] 
\small {\sc $^1$Oslo Centre for Biostatistics and Epidemiology (OCBE) \\
   and $^2$Department of Mathematics, University of Oslo} \\[0.3em]
\small {\sc {\heute}}\par
\end{centering}
\endgroup


\begin{abstract}
  \small{In various image analysis contexts there is a need to accurately
    estimate the background variance. A stylised version of that problem,
    perhaps after pre-processing of images, 
    is to estimate the standard deviation $\sigma$ for the $\eps_i$
    in the model $y_i=m(x_i)+\eps_i$ for $i=1,\ldots,n$,
    with $m(x)$ the signal, assumed to be locally constant, along with
    jumps and fluctuations. Inside the software package
    AQuA (Astrocyte Quantitative Analysis),
    developed in \citet*{Wangetal19} for analysing
    neurotransmitter fluorescence dynamics for single-cell
    and population-level physiology, the estimator 
    $\hatt\sigma_W=M_n/c$ is constructed,
    in terms of the median $M_n$ of absolute differences $w_i=|y_{i+1}-y_i|$,
    with $c=\sqrt{2}$.
    This is part of their crucial Step 1, the first of 8 steps,
    though without any theoretical examination of its behaviour.
    In this note we properly sort out its distribution,
    first establishing that one should rather use the correction factor
    $c=0.9539$ for consistency of the normal $\sigma$.
    This better estimator is employed in other software,
    developed for similar and specifically for Ca$^{2+}$
    signal image analysis process purposes, in
    \citet{Bjoernstad21}, \citet{Aabjoern22}, \citet{Bojarskaite23}. 
    Applying limit distribution theory for quantile processes
    we propose the improved practical estimator
    $Q_{0.75}/1.629$, using the 0.75 quantile of the $w_i$. 
    
}
\end{abstract}

\smallskip 
\begin{center}
\fbox{%
\parbox{0.77\textwidth}{%
    The authors have not used any AI for this work, neither for
    the mathematics, statistics, probability,
    or generating of ideas, nor, in particular, 
    for the actual writing of words, sentences, paragraphs. 
    Vivat Homo sapiens. 
  } %
}
\end{center} 

\section{Introduction}
\label{section:intro} 

Suppose $y_1,\ldots,y_n$ are observed, with a constant
variability but with a mean level function assumed 
to have flat parts, along with jumps and fluctuations. 
How should one then estimate the variability level?
Such problems occur in a wide variety of statistical
application domains, in particular for analyses of
images, perhaps monitored over time. 


\begin{figure}[h]
\centering
\includegraphics[scale=0.40]{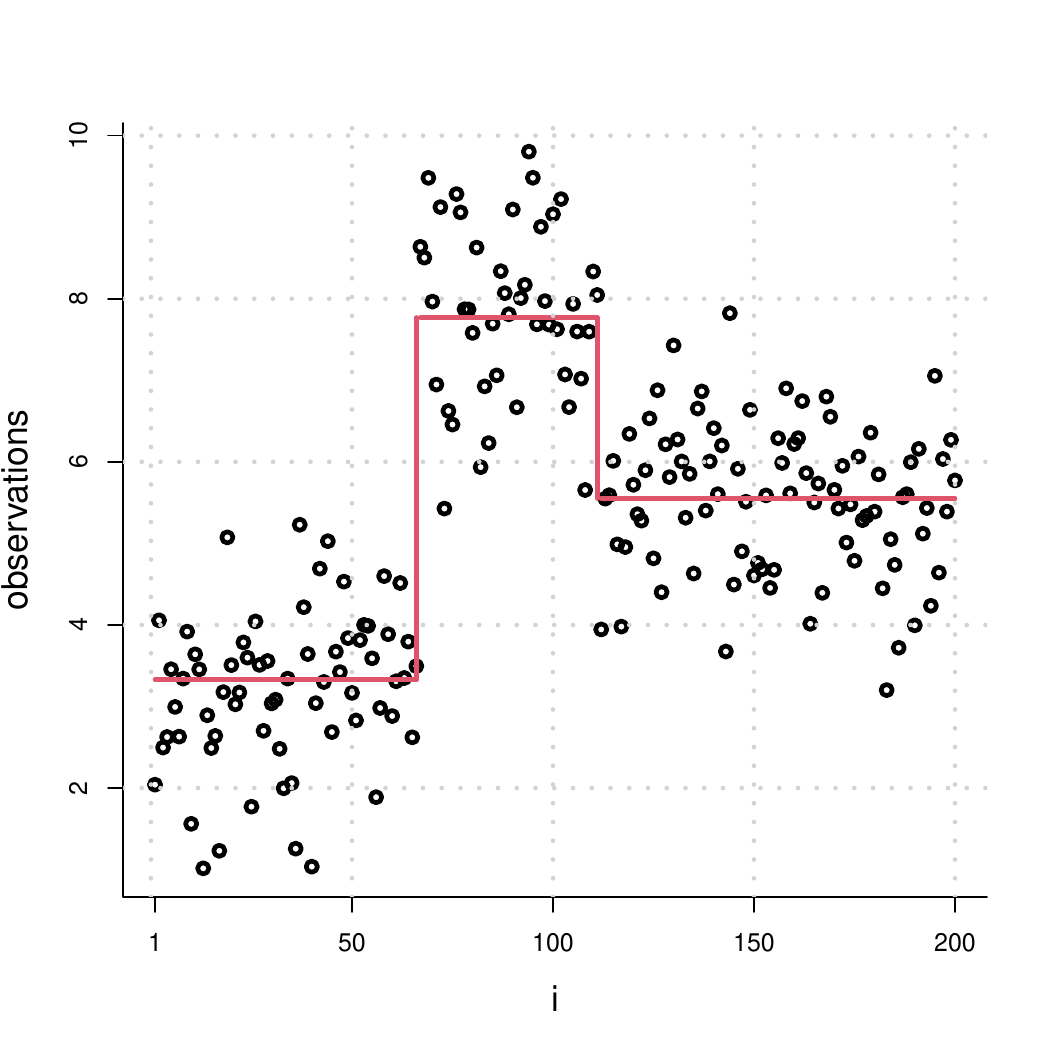}
\includegraphics[scale=0.40]{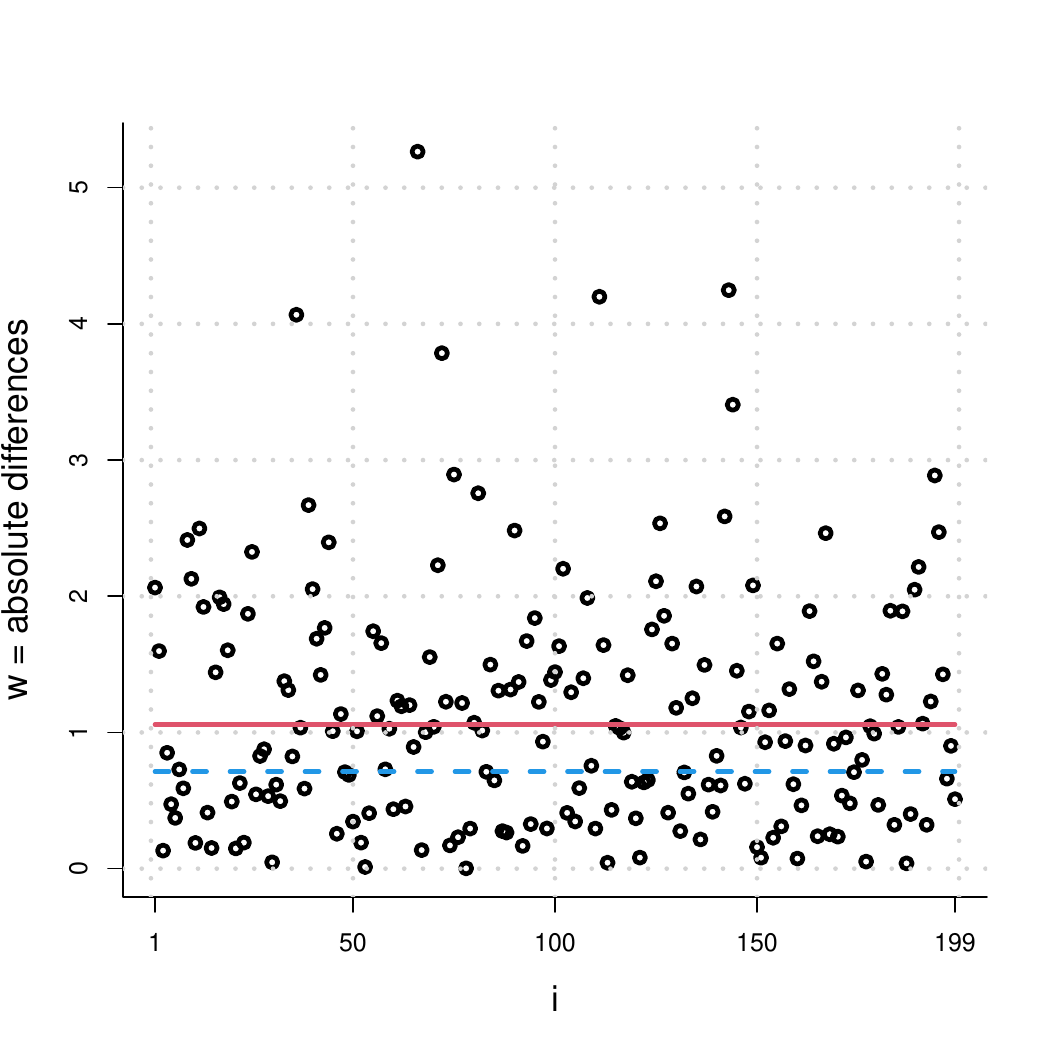}
\caption{Simulated observations exhibiting jumps in mean level,
   with mean function plotted in red, 
   but with a constant variability (left panel). 
   The right panel displays the absolute differences
   $w_i=|y_{i+1}-y_i|$, with median along with two horizontal lines;
   the full red line is for $\hatt\sigma=M_n/0.9539$,
   close to the true $\sigma=1$, the dashed blue line
   is for the Wang et al.~estimator $\hatt\sigma_W=M_n/\sqrt{2}$,
   which is far too low.}
\label{figure:fig12}
\end{figure}

A formalised setup might start with 
\beq
y_i=m(x_i)+\eps_i, \quad {\rm for\ }i=1,\ldots,n
\label{eq:model}
\eeq
where the $\eps_i$ are i.i.d.~$\N(0,\sigma^2)$, 
along with relevant conditions on the signal function $m(x)$,
as alluded to above. 
It is clear that the traditional estimator,
the empirical standard deviation, risks serious
overshooting of the real $\sigma$. Not surprisingly
several types of alternative estimators can 
be constructed, depending also on the level
and type of anticipated jumps. 
Relevant estimators, essentially able to give 
the right level of variability in the presence of 
such jumps in level, include those 
proposed and studied in \citet*{Neumannetal41}, 
\citet{RousseeuwCroux93}. 

The present paper is concerned with an estimator
proposed by \citet{Wangetal19}, in a study involving 
many such time series with jumps and fluctuations.
These authors have constructed a full 
Astrocyte Quantitative Analysis (AQuA) software package
for analysing complex calcium and neurotransmitter activity
in fluorescence imagining datasets, both {\it ex vivo}
and {\it in vivo}. As part of their Stage 1 processing
of such data, they employ the following robust scale estimator, 
for each of many parallel observation sequences; 
with $y_1,\ldots,y_n$, the algorithm computes 
\beq
\hatt\sigma=\med\{w_1,\ldots,w_{n-1}\}/c
\quad {\rm with\ }
w_i=|y_{i+1}-y_i|,  
\label{eq:wangestimator} 
\eeq  
and $c$ an appropriate scaling factor. \citet{Wangetal19} 
use the factor $c=\sqrt{2}$. The paper mentioned does not
give any properties of their estimator, which we will tend
to below, however. In particular, we rather prefer 
the scaling factor needed to make the estimator 
consistent for $\sigma$ under normal conditions; 
it turns out that $c=\sqrt{2}\,\Gamma_1^{-1}(\half)^{1/2}=0.9539$ 
is then required. 

In this note  we study the estimator (\ref{eq:wangestimator})
more carefully, and in particular exhibit the precise 
limit distribution of $\rootn(\hatt\sigma-\sigma)$
under normal conditions. This enables inferences
to be made, such as confidence distributions,
comparisons of different sequences, etc. 
Reaching such results requires an apparatus of 
empirical process, along with quantile versions
of these. These techniques, with results, are 
given in Section \ref{section:mainstuff}. 
In particular it is found that the (\ref{eq:wangestimator})
estimator is rather inefficient, under flat model conditions,
leading to confidence intervals 1.85 times wider
than with the usual empirical standard deviation. 

In Section \ref{section:pquantestimator} we
propose and analyse a wider class of estimators,
related to taking a general $p$ quantile of the 
$w_1,\ldots,w_{n-1}$ absolute differences, rather than 
$p=0.50$. It turns out, via analysis of the 
corresponding $\rootn(\hatt\sigma_p-\sigma)$
distribution, that its variance goes down with increasing $p$,
to its best value of $p=0.865$, and then
turning up after that optimal position. 
A few concluding comments are placed in Section \ref{section:concluding}. 

We note that the better variance estimator,
improving on the Wang et al.~estimator, 
constructed via theory developed below,
is employed in software, developed for various 
image analysis process purposes, in
\citet{Bjoernstad21}, \citet{Aabjoern22}, \citet{Bojarskaite23}.
In particular, pages 4-5 in the Begonia paper \citet{Bjoernstad21}
outline the specific role the $\sigma$ estimator
plays in ensuing statistical image analysis
for monitoring the intact brains of awake mice without
the dampening effects of anesthesia, revealing
a rich repertoire of astrocytic Ca$^{2+}$ signals. 


\def\emp{{\rm emp}}

\section{Limit distribution for the Wang et al.~estimator} 
\label{section:mainstuff} 

Let $y_1,\ldots,y_n$ be i.i.d.~$\N(\xi,\sigma^2)$.
We shall investigate the limit distribution 
for the $\hatt\sigma$ of (\ref{eq:wangestimator}). 
We note first, for comparison purposes, that 
the standard method, with the empirical 
standard deviation $\hatt\sigma_\emp$, has 
\beq
\rootn(\hatt\sigma_\emp^2-\sigma^2)\arr_d \N(0,2\sigma^4)
\quad {\rm and\ hence} \quad 
\rootn(\hatt\sigma_\emp-\sigma)\arr_d \N(0,\half\sigma^2). 
\label{eq:standardestimator}
\eeq 
The first result is classic, following from the 
$\hatt\sigma^2/\sigma^2\sim\chi^2_{n-1}$ property,
and the second follows by applying the delta method.
One question is how much bigger the limiting variance will be,
for $\rootn(\hatt\sigma-\sigma)$,
compared to the best-under-model $\half\sigma^2$.  

For the successive absolute differences $w_i=|y_{i+1}-y_i|$  
of (\ref{eq:wangestimator}), we have $w_i^2 \sim 2\sigma^2\chi^2_1$
under model conditions, so the median of a high number of such, 
whether independent or not, tends to
$\sqrt{2}\sigma\,\Gamma_1^{-1}(\half)^{1/2}=0.9539\,\sigma$.  
Here we write $\Gamma_1(\cdot)$ for the cumulative distribution function 
(cdf) of the $\chi^2_1$, with density $\gamma_1(\cdot)$. 
The correct scaling factor to use for 
\beq
\hatt\sigma=\med(w_1/c,\ldots,w_{n-1}/c) 
\label{eq:sigmahat}
\eeq
is hence $c=0.9539$. 
We shall show below that under normal i.i.d.~conditions, 
\beq
\rootn(\hatt\sigma-\sigma)
   \arr_d\N(0,\sigma^2\tau^2),
\quad {\rm where\ }
\tau={ (\quart+2v_1)^{1/2}\over \gamma_1(\half c^2)c^2},   
\label{eq:mainthingsection2}
\eeq 
involving 
\beqn
v_1 = \Pr\{|z_2-z_1|\le c,|z_3-z_2|\le c \}-\quart
    = \Pr(|N_1|\le c/\sqrt{2},|N_2|\le c/\sqrt{2})-\quart, 
\eeqn 
for independent standard normals $z_1,z_2,z_3$, and where the second
formula writes this in terms of standard normals $N_1,N_2$
with correlation $-\half$. 
The numerical value of $\tau$ is 1.285 which 
is as much as 1.816 times bigger than the $1/\sqrt{2}=0.7071$
attained by the optimal $\hatt\sigma_\emp$ of 
(\ref{eq:standardestimator}). Better estimators,
generalising $\hatt\sigma$ of (\ref{eq:sigmahat}), 
with less loss of efficiency at the standard 
flat normal model, will be worked with in 
Section \ref{section:pquantestimator}.  

It suffices to study the case of $\sigma=1$,
and then to multiply the resulting limit distribution
appropriately. 
Let $M_n=\med(w_1,\ldots,w_{n-1})$. 
Consider the empirical cdf $G_n(t)$ for the $w_i$,
converging in probability to
$G(t)=\Pr(\sqrt{2}|N(0,1)|\le t)=\Gamma_1(\half t^2)$,
with density $g(t)=\gamma_1(\half t^2)t$. 
The plan is to first establish a process convergence result of the type 
\beq
Z_n(t)=\rootn\{G_n(t)-G(t)\}\arr_d Z(t),
\label{eq:notdonsker}
\eeq 
where the limit is a zero-mean normal process 
with a certain covariance function. 
It then follows, via general results of \citet{DossGill92}, 
that 
\beq
\rootn\{G_n^{-1}(p)-G^{-1}(p)\}
   \arr_d -{Z(G^{-1}(p))\over g(G^{-1}(p))}, 
\label{eq:quantilelimit}
\eeq 
for the full associated quantile process. In particular,
with $c=G^{-1}(\half)=0.9539$, 
we will then have found 
\beqn
\rootn(M_n/c-1)\arr_d -Z(c)/\{g(c)c\},
\quad {\rm where\ }
g(c)=\gamma_1(\half c^2)c=0.4494. 
\eeqn 
This is a normal $(0,\sigma^2\tau^2)$ limit, with
$\tau^2=\{\Var\,Z(c)\}/ \{g(c)^2c^2\}$, with the variance of $Z$
found below. This is what leads to (\ref{eq:mainthingsection2}). 


It remains to establish (\ref{eq:notdonsker}),
which is not covered by the classic Donsker theorem,
since the $w_i$ variables are not independent. 
The result is however secured via the general 
methods and results of \citet[Ch.~4]{Billingsley68};
see in particular his Section 22 and Theorem 22.1.
Relevant further material is in 
\citet[Ch.~13]{HjortStoltenberg27} for functional
central limit theorems for 1-dependence processes. 
There is a well-defined limit process $Z(\cdot)$,
and its covariation function is not difficult 
to write down. In our situation we are primarily
concerned with the consequence (\ref{eq:quantilelimit})
of (\ref{eq:notdonsker}), so we need only the variance. 
Indeed, 
\beq
Z_n(t)\arr_d Z(t)\sim\N(0,v_0(t)+2v_1(t)), 
\label{eq:ZntoZ} 
\eeq 
where 
$v_0(t)=\Var\,I(w_i\le t)=G(t)\{1-G(t)\}$
and 
\beqn
v_1(t)
&=&\cov\,\{I(w_i\le t),I(w_{i+1}\le t)\} \\
&=&\Pr(w_i\le t,w_{i+1}\le t)-G(t)^2 \\
&=& \Pr(|N_1|\le t/\sqrt{2},|N_2|\le t/\sqrt{2}), 
\eeqn 
in terms of standard normals $N_1$ and $N_2$
with correlation $-\half$, amenable to numerical integration
or simulation for its computation.
Numerical work now verifies that (\ref{eq:mainthingsection2})
holds, for the estimator $\med(w_1,\ldots,w_n)/c$, with $\tau=1.285$. 


\section{A more general quantile based estimator} 
\label{section:pquantestimator} 

The $\hatt\sigma$ estimator of (\ref{eq:sigmahat}) 
uses the median of the observed $w_i=|y_{i+1}-y_i|$, 
suitably scaled. We have seen that there is 
a significant loss of efficiency, under flat model conditions.
A generalisation of the basic  estimator is 
\beq
\hatt\sigma_p=\quantile_p (w_1/c_p,\ldots,w_{n-1}/c_p)
   =G_n^{-1}(p)/c_p, 
\label{eq:generalsigmahat}
\eeq
the empirical $p$-quantile of the observed $w_i/c_p$,
for a given quantile value $p$. This requires again 
an appropriate scaling factor $c_p$ to achieve consistency,
$\hatt\sigma_p\arr_\pr\sigma$. 
In fact, from $G(t)=\Gamma_1(\half t^2)$ above, set equal to $p$,
we need $c_p=\sqrt{2}\,\Gamma_1^{-1}(p)^{1/2}$ for consistency. 
From efforts above, mutatis mutandis, we learn
for the estimator $\hatt\sigma_p=G_n^{-1}(p)/c_p$ that 
\beqn
\rootn(\hatt\sigma_p-\sigma)\arr_d
   -\sigma {Z(c_p)\over g(c_p)c_p} \sim\N(0,\sigma^2\tau_p^2),
\eeqn 
with
\beq
\tau_p={1\over g(c_p)c_p} \{v_0(c_p)+2v_1(c_p)\}^{1/2}
   = { \{p(1-p)+2(B_p-p^2)\}^{1/2} \over \gamma_1(\Gamma_1^{-1}(p)) c_p^2},  
\label{eq:taup}
\eeq
in terms of $B_p=\Pr(|N_1|\le c_p/\sqrt{2},|N_2|\le c_p/\sqrt{2})$,
with the two standard normals having correlation $-\half$. 

\begin{figure}[h]
\centering
\includegraphics[scale=0.6]{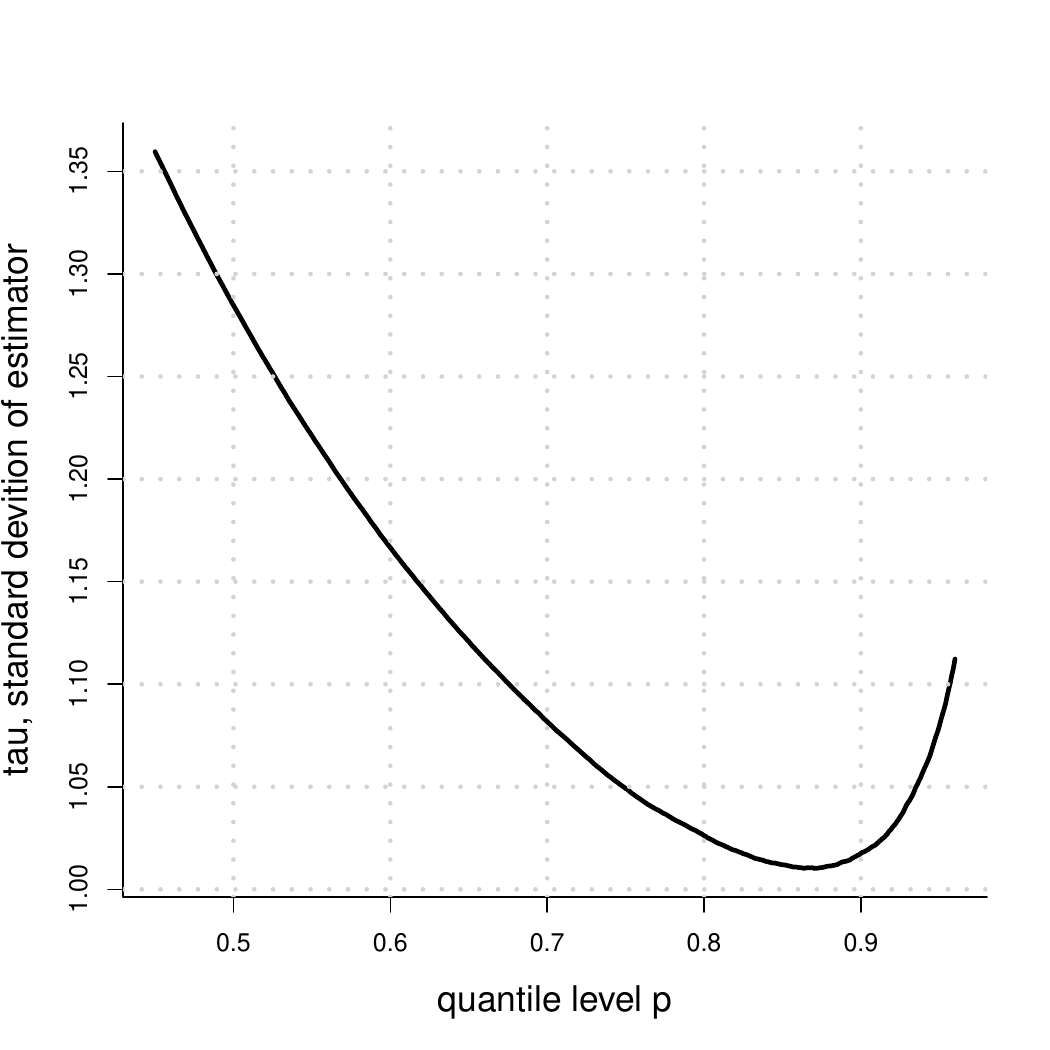}
\caption{The limit distribution standard deviation $\tau_p$
   of (\ref{eq:taup}), as a function of quantile level $p$,
   for the estimator $\hatt\sigma_p$ of (\ref{eq:generalsigmahat}).
   The minimum occurs of $p=0.86$. The \citet{Wangetal19} 
   estimator corresponds to $p=0.50$, where the standard deviation
   is 1.273 times bigger.}
\label{figure:fig11}
\end{figure}

The $\tau_p$ is displayed, as a function of quantile level $p$,
in Figure \ref{figure:fig11}. We learn that it pays off 
in efficiency to use a higher value than $p=0.50$.
The minimum occurs for $p=0.086$, 
with a minimum $\tau_p=1.009$. This ought to be compared
with $1/\sqrt{2}=0.7071$, which is the smallest possible
value when estimating $\sigma$ under normal i.i.d.~conditions.
The strength of the $\hatt\sigma_p$ estimator of 
(\ref{eq:generalsigmahat}) is of course that it is
robust against jumps and fluctuations in the underlying signal function.

For practical applications we propose computing the 0.75 quantile
$Q_{0.75}=G_n^{-1}(0.75)$ of the $w_i=|y_{i+1}-y_i|$,
and then using the estimator
$\hatt\sigma_{0.75}=Q_{0.75}/c_{0.75}=Q_{0.75}/1.629$.
We have $\rootn(\hatt\sigma_{0.75}/\sigma-1)\arr_d\N(0,1.049^2)$,
from which confidence intervals, and in fact a full confidence curve,
can be found. 


\section{Robustness of the robust estimator against a few sudden jumps} 
\label{section:concluding}

Large-sample formulae have been derived above, with the key result
being $Z_n=\rootn(G_n-G)\arr_d Z$ for the absolute differences,
under normal i.i.d.~conditions for the data. 
If there are sudden mean shifts, this may affect the limits.
Writing $\mu_1,\ldots,\mu_n$ for the mean, we have
\beqn
w_i=|y_{i+1}-y_i|=|\mu_{i+1}-\mu_i+\eps_{i+1}-\eps_i|, 
\eeqn 
with the $\eps_i$ being i.i.d.~normal. We learn from this that as long as
there are $s_n$ mean jumps, among the $n$, and $s_n/\rootn\arr0$,
then the $Z_n\arr_d Z$ result stands unaffected, hence also its
corollaries. 

Results above rest on the normality assumption, specifically
for the variance in (\ref{eq:ZntoZ}).
If $y_i=m(x_i)+\sigma\eps_i$, with the $\eps_i$ i.i.d.~with zero mean
and unit variance, the limit variance takes the form
\beqn
\Var\,Z(t)=G(t)\{1-G(t)\}
   + 2 \{\Pr(|\eps_2-\eps_1|\le t,|\eps_3-\eps_2|\le t)-G(t)^2\}, 
\eeqn
with $G(t)=\Pr(|\eps_2-\eps_1|\le t)$,
specialising to the above for the normal case. 
Brief inspection indicates that the limit distribution for
$\rootn(\hatt\sigma_{0.75}-\sigma)$ is not much affected
by non-normality. 

The basis for being to derive these results
is (\ref{eq:notdonsker}); then the quantile results
follow, as shown via the techniques of \citet{DossGill92}.
We may also obtain results under non-i.i.d.~assumptions
for the $y_i$ observations, as with time series
models with jumps in them. In such situations
the theory will go through, but with more complicated
variance formulae for the $Z_n\arr_d Z$. 


%

%

\bibliographystyle{biometrika}
\bibliography{diverse_bibliography2026}

\end{document}